\documentclass[10pt,twoside,a4paper]{amsart}
\usepackage[english]{babel}
\usepackage[utf8]{inputenc}

\usepackage{amsfonts,amsthm,amssymb,amsmath,amsthm,latexsym,wasysym,stmaryrd,mathrsfs,dsfont,txfonts}
\usepackage{accents,multirow,rotating,subfigure,graphicx,bigstrut,lscape,multicol,fancyhdr,enumerate,accents,mathtools,exscale}

\usepackage{pdftricks}
\usepackage[pdftex,all]{xy}

\usepackage{hyperref}
\usepackage{appendix}
\usepackage[justification=centering]{caption}

\theoremstyle{theorem}
\newtheorem {thm}{Theorem}[section]
\newtheorem*{thm*}{Theorem}
\newtheorem*{thmdef}{Theorem--Definition}
\newtheorem*{thmdef*}{Theorem--Definition}
\newtheorem {lemma}[thm]{Lemma}
\newtheorem*{lemma*}{Lemma}

\newtheorem*{prop*}{Proposition}

\newtheorem*{cor*}{Corollary}

\newtheorem*{conjecture*}{Conjecture}
\theoremstyle{definition}
\newtheorem {defi}[thm]{Definition}
\newtheorem*{defi*}{Definition}

\newtheorem*{nota*}{Notation}
\theoremstyle{remark}
\newtheorem {remark}[thm]{Remark}
\newtheorem*{remark*}{Remark}

\newtheorem*{warning*}{Warning}

\newtheorem*{remarques*}{Remarks}

\newtheorem*{warnings*}{Warnings}

\newtheorem*{convention*}{Convention}
\newtheorem {example}[thm]{Example}

\newtheorem*{question*}{Question}

\newtheorem*{questions*}{Questions}

\newtheorem*{fact*}{Fact}
\newtheorem*{acknowledgments}{Acknowledgments}

\def\omu{\overline{\mu}}

\def\Z{{\mathds Z}}

\def\lk{\textrm{lk}}
\def\F{\textrm{F}}

\newcommand{\dessin}[2]{
  \vcenter{\hbox{\includegraphics[height=#1]{#2.pdf}}}}

\begin{document}

\title{Linking number and Milnor invariants}
\author[J.B. Meilhan]{Jean-Baptiste Meilhan}
         \address{Universit\'e Grenoble Alpes, IF, 38000 Grenoble, France}
         \email{jean-baptiste.meilhan@univ-grenoble-alpes.fr}

\begin{abstract}
This is a concise overview of the definitions and properties of the linking number and its higher-order generalization, Milnor invariants. 
\end{abstract}

\maketitle

\section{Introduction}
The linking number is probably the oldest invariant of knot theory. 
In 1833, several decades before the seminal works of physicists Tait and Thompson on knot tabulation, 
the work of Gauss on electrodynamics led him to formulate the number of ``intertwinings of two closed or endless curves'' as
 \begin{equation}\label{gauss_lk}
 \frac{1}{4\pi} \int \int \frac{(x'-x)(dydz'-dzdy') + (y'-y)(dzdx'-dxdz')+(z'-z)(dxdy'-dydx')}{[(x'-x)^2 + (y'-y)^2 + (z'-z)^2]^{3/2}}, 
\end{equation}
where $x, y, z$ (resp. $x', y', z'$) are coordinates on the first (resp. second) curve. 
This number, now called 
\emph{linking number}, is an invariant of $2$--component links,\footnote{In this note, all links will be ordered and oriented tame links in the $3$-sphere.  }
which typically distinguishes the trivial link and the Hopf link, shown below. 
\[ \includegraphics[scale=0.6]{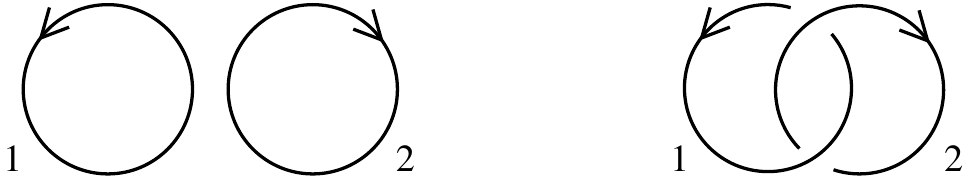} \]
Alternatively, the linking number can be simply defined in terms of the first homology group of the link complement (see Thm.~\ref{def:lk1}). 
But, rather than the abelianization of the fundamental group, one can consider finer quotients, for instance by the successive terms of its lower central series, to construct more subtle link invariants. 
This is the basic idea upon which Milnor based his work on higher order linking numbers, now called \emph{Milnor $\overline{\mu}$-invariants}. 

We review in Section  \ref{sec:lk} several definitions and key properties of the linking number. 
In Section \ref{sec:mu}, we shall give a precise definition of Milnor invariants and explore some of their properties, generalizing those of the linking number. 
In the final Section \ref{sec:more}, we briefly review further known results on Milnor invariants. 

\section{The Linking Number} \label{sec:lk}

The linking number has been studied from multiple angles and has thus been given many equivalent definitions of various natures. 
We already saw in (\ref{gauss_lk}) the original definition due to Gauss; 
let us review below a few others.  
\subsection{Definitions}
For the rest of this section, let $L=L_1\cup L_2$ be a $2$ component link, and let  $\lk(L)$ denote the linking number of $L$.

\vspace{.1cm}
\hspace{-.53cm}
\parbox[l]{14cm}{
$\quad$ Recall that the first homology group of the complement $S^3\setminus L_1$ of $L_1$ is the infinite cyclic  $\quad$ 
group generated by the class $[m_1]$ 
of a meridian $m_1$, as shown on the right. 
\begin{thm}\label{def:lk1}
The linking number of $L$ is the integer $k$ such that $[L_2]=k
[m_1]\in H_1(S^3\setminus L_1;\Z)$.  
\end{thm}
   } 
  \parbox[r]{0.7cm}{
  }
  \parbox[r]{1cm}{
    \vspace{-.4cm}
  \includegraphics[scale=0.8]{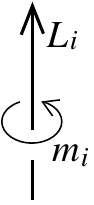}
  }

Recall now that any knot bounds orientable surfaces in $S^3$, called Seifert surfaces. 
Now, given a Seifert surface $S$ for $L_2$, one may assume up to isotopy that $L_1$ intersects $S$ at finitely many transverse points.  At each intersection point, the orientation of $L_1$ and that of $S$ form an oriented $3$-frame, hence a sign using the right-hand rule; the \emph{algebraic intersection} of $L_1$ and $S$ is the sum of these signs over all intersection points. 
\begin{thm}\label{def:lk2}
The linking number of $L$ is the algebraic intersection of $L_1$ with a Seifert surface for $L_2$.  
\end{thm}

Given a regular diagram of $L$, the linking number is simply given by the number of crossings where $L_1$ passes over $L_2$, counted with signs as follow.  
\begin{thm}\label{def:lk3}
\[ 
\lk(L) = \left( \textrm{number of }
\begin{array}{l}
\includegraphics{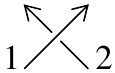}
\end{array} \right) 
 - \left( \textrm{ number of }
\begin{array}{l}
\includegraphics{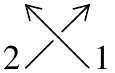}. 
\end{array} \right)
 \]
\end{thm}

The proof of the above results, \emph{i.e.} the fact that these definitions are all equivalent and equivalent to (\ref{gauss_lk}), can be found in \cite[\S.~5.D]{Rolfsen}, along with several further definitions. See also \cite{Ricca} for details. 
\begin{remark}
As part of Theorems \ref{def:lk2} and \ref{def:lk3}, the linking number does not depend on the choices involved in these definitions, namely the choice of a Seifert surface for $L_2$ and of a diagram of $L$, respectively. 
\end{remark}

\subsection{Basic properties}

It is rather clear, using any of the above definitions, that reversing the orientation of either component of $L$ changes the sign of $\lk(L)$. 

Also clear from Gauss formula (\ref{gauss_lk}) is the fact 
that the linking number is symmetric : 
\begin{equation}\label{sym_lk}
  \lk(L_1,L_2)=\lk(L_2,L_1). 
\end{equation}
\noindent This can also be seen, for example, from the diagrammatic definition (Thm.~\ref{def:lk3}), 
by considering two projections, on two parallel planes which are `on either sides' of $L$, so that crossings where $1$ overpasses  $2$ in one diagram are in one-to-one correspondence with crossings where $2$ overpasses $1$ in the other, with same sign. 

The symmetry property (\ref{sym_lk}) allows for a symmetrized version of Theorem~\ref{def:lk3}: 
\[ \lk(L) = \frac{1}{2}\left(\big( \textrm{number of }
\begin{array}{l}
\includegraphics{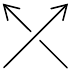}
\end{array} \big) 
 - \big( \textrm{ number of }
\begin{array}{l}
\includegraphics{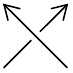}  
\end{array} \big) \right), 
 \]
\noindent where the sum runs over \emph{all} crossings involving a strand of component $1$ and a strand of component $2$. 

\subsection{Two classification results} \label{sec:classiflk}

\begin{figure}[h]
$ \dessin{1.5cm}{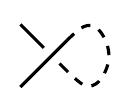} \longleftrightarrow
      \dessin{1.5cm}{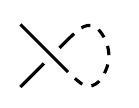}
\qquad 
  \dessin{1.5cm}{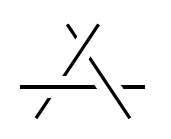} \longleftrightarrow
        \dessin{1.5cm}{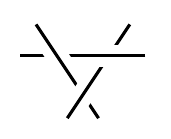}$
  \caption{A self-crossing change (left) and a $\Delta$-move (right)}
  \label{fig:local}
\end{figure}

The following notion was introduced by Milnor in the fifties \cite{Milnor}.
\begin{defi}\label{def:lh}
 Two links are \emph{link-homotopic} if they are related by a sequence of ambient isotopies and self-crossing changes, \emph{i.e.} crossing changes involving two strands of a same component (see Figure \ref{fig:local}).
\end{defi}
The idea behind this notion is that, as a first approximation to the general study of links, working up to link-homotopy allows to unknot each individual component of a link, and only records their `mutual interractions' -- this is in a sense studying `linking modulo knotting'. 

Observe that the linking number is invariant under link-homotopy. 
Furthermore, it is not too difficult to check that any $2$--component link is link-homotopic to an iterrated Hopf link, \emph{i.e.}  the closure of the pure braid $\sigma_1^{2n}$ for some $n\in \Z$. This number $n$ is precisely the linking number of the original link, thus showing: 
\begin{thm}\label{prop:lh}
 The linking number classifies $2$-component links up to link-homotopy. 
\end{thm}
Generalizing this result to a higher number of components, however, requires additional invariants, and this was one of the driving motivations that led Milnor to develop his $\omu$-invariants, see Section \ref{sec:lhc}. 

The equivalence relation on links with an arbitrary number of components which is classified by the linking number is the 
\emph{$\Delta$-equivalence}, which is generated by the $\Delta$-move, shown on the right-hand side of  Figure \ref{fig:local}. 
Indeed, the following is due H.~Murakami and Y.~Nakanishi \cite{MN}, see also \cite{Matveev}. 
\begin{thm}\label{thm:delta}
 Two links $\cup_{i=1}^n L_i$  and $\cup_{i=1}^n L'_i$ ($n\ge 2$) are 
 $\Delta$-equivalent if, and only if $\lk(L_i,L_j)=\lk(L'_i,L'_j)$ for all $i,j$ such that $1\le i<j\le n$.
\end{thm}
\noindent Notice in particular that the $\Delta$-move is an unknotting operation, meaning that any knot can be made trivial by a finite sequence of isotopies and $\Delta$-moves. 

\section{Milnor invariants} \label{sec:mu}

In this whole section, let $L=L_1\cup \cdots \cup L_n$ be an $n$-component link, for some fixed $n\ge 2$. 

\subsection{Definition of Milnor $\omu$-invariants}
In his master's and doctoral theses,  supervised by R.~Fox, Milnor defined numerical invariants extracted from the peripheral system of a link, which widely generalize the linking number. 

Denote by $X$ the complement of an open tubular neighborhood of $L$.
Pick a point in the interior of $X$, and denote by $\pi$ the fundamental group of $X$ based at this point.
It is well-known that, given a diagram of $L$, we can write an explicit presentation of $\pi$, called Wirtinger presentation, where each arc in the diagram provides a generator, and each crossing yields a relation. 
Such presentations, however, are in general difficult to work with, owing to their 
large number of generators. 
We introduce below a family of quotients of $\pi$ for which a much simpler presentation can be given, and which still retain rich topological information on the link $L$. 

The \emph{lower central series} $(\Gamma_k G)_{k\ge 1}$ of a group $G$ is the nested family of subgroups defined inductively by 
 $$ \Gamma_1 G = G\quad \textrm{and}\quad  \Gamma_{k+1} G = [G,\Gamma_k G]. $$
The \emph{$k$th nilpotent quotient} of the group $G$ is the quotient $G / \Gamma_k G$. 

Let us fix, from now on, a value of $k\ge 2$. 
For each $i\in\{1,\cdots,n\}$, consider two elements $m_i,\lambda_i\in \pi / \Gamma_k \pi$, where $m_i$ represents a \emph{choice} of a meridian for $L_i$ in $X$,  
and where $\lambda_i$ is a word representing the $i$th \emph{preferred longitude} of $L_i$,  \emph{i.e.}  a parallel copy of $L_i$ in $X$ having linking number zero with $L_i$. 
Denote also by $F$ the free group $F(m_1,\cdots,m_n)$.
\begin{thm}[Chen-Milnor Theorem]\label{thm:CM}
The $k$th nilpotent quotient of $\pi$ has a presentation given by 
 $$ \pi / \Gamma_k \pi = \left\langle m_1,\ldots, m_n\, \vert\,  m_i\lambda_im_i^{-1}\lambda_i^{-1}\, (i=1,\cdots,n)\, ; \, \Gamma_k\F \right\rangle. $$
\end{thm}

In order to extract numerical invariants from the nilpotent quotients of $\pi$, we consider the \emph{Magnus expansion} \cite{MKS} of each $\lambda_j$, which is an element $E(\lambda_j)$ of $\mathbb{Z}\langle \langle X_1, . . . ,X_n\rangle \rangle$, the ring of formal power series in non commuting variables $X_1, \ldots, X_n$, obtained by the substitution  
\[ m_i \mapsto  1 + X_i \, \, \, \textrm{ and }\, \, \, 
 m_i^{-1}  \mapsto  1 - X_i + X_i^2 - X_i^3 + \ldots 
\]

We denote by $\mu_L(i_1 i_2 \ldots i_p j)$ the coefficient of $X_{i_1} X_{i_2}\ldots X_{i_p}$ in the Magnus expansion of $\lambda_j$:
\begin{equation}\label{eq:E}
 E(\lambda_j) = 1 + \sum_{i_1 i_2 \ldots i_p} \mu_L(i_1 i_2 \ldots i_p j) X_{i_1} X_{i_2}\ldots X_{i_p}
\end{equation}
This coefficient is \emph{not}, in general, an invariant of the link, as it depends upon the choices made in this construction (Essentially, 
picking a system of based meridians for $L$).  We can, however, promote these number to genuine invariants by regarding them modulo the following indeterminacy.  
\begin{defi}
Given a sequence $I$ of indices in $\{1, ... , n\}$, let $\Delta_L(I)$ be the greatest common divisor of the coefficients $\mu_L(I')$, for all sequences  $I'$ obtained from $I$ by deleting at least one index and permuting cyclicly.
\end{defi}
For example, $\Delta_L(123) = \textrm{gcd}\left\{ \mu_L(12), \mu_L(21), \mu_L(13), \mu_L(31), \mu_L(23), \mu_L(32) \right\}$.
\begin{thmdef}[Milnor \cite{Milnor2}]
The residue class
 $$ \overline{\mu}_L(I) \equiv \mu_L(I)\textrm{ mod }\Delta_L(I) $$
is an invariant of ambient isotopy of $L$, for any sequence $I$ of $\le k$ integers, 
called a \emph{Milnor invariant} of $L$. The number of indices in $I$ is called the \emph{length}. 
\end{thmdef}
Before giving some examples, in Section \ref{sec:exmu}, a few remarks are in order. 
\begin{remark}\label{rem:length1}
Observe from (\ref{eq:E}) that $\mu_L(I)$ is only defined for a sequence $I$ of \emph{two} or more indices.  We set $\mu_L(i)=0$ for all $i$ as a convention. 
\end{remark}
\begin{remark}\label{rem:firstnonzero}
If all Milnor invariants of $L$ of length $<k$ are zero, then those of length $k$ are well defined integers. These \emph{first non-vanishing Milnor invariants} are thus much easier to study in practice and, as a matter of fact, a large proportion of the literature on  Milnor link invariants focusses on them. 
\end{remark}
\begin{remark}\label{rem:k}
Working with the $k$th nilpotent quotient of $\pi$ only allows to define Milnor invariants up to length $k$. But this $k$ can be chosen arbitrarily large, so this is no restriction. 
\end{remark}

\subsection{First examples}\label{sec:exmu}

By Remark \ref{rem:length1}, length $2$ Milnor invariants are well-defined over $\Z$. 
In order to define them, it suffices to work in the second nilpotent quotient, \emph{i.e.} in $\pi/\Gamma_2\pi=H_1(X;\Z)$. But the $j$th preferred longitude in $H_1(X;\Z)$ is given by  
$\lambda_j  = \sum_{i\neq j} \lk(L_i,L_j) m_i$, 
and we have  
\[ E(\lambda_j) = 1+ \sum_{i\neq j} \lk(L_i,L_j) X_i + \textrm{terms of degree $\ge 2$}. \]
Hence length $2$ Milnor invariants are exactly the pairwise linking numbers of a link. This justifies regarding length $\ge 3$  Milnor invariants as `higher order linking numbers'. 
In order to get a grasp on these, let us consider an elementary example. 

\begin{example} \label{ex:bor}
Consider the Boromean rings $B$ , as illustrated below.\footnote{
There, we use the usual convention that a small arrow underpassing an arc of the diagram represents a loop, based at the reader's eye and going straight down to the projection plane, enlacing positively the arc and going straight back to the basepoint. } 
The Wirtinger presentation for $\pi$ is    

\vspace{.15cm}
\hspace{-.43cm}
\parbox[l]{10.7cm}{
$\quad \pi = \left\langle \begin{array}{c|ccc}
m_1\, ,\, m_2\,,\, m_3\, & m_3m_1m_3^{-1}n_1^{-1} \, \,  ;&   m_1m_2m_1^{-1}n_2^{-1} \, \, ;& m_2m_3m_2^{-1}n_3^{-1}\\
n_1\, ,\, n_2\,,\, n_3\, &   n_1^{-1}n_2n_1m_2^{-1}\, \,  ;&   n_2^{-1}n_3n_2m_3^{-1}   \, \, ;&
 n_3^{-1}n_1n_3m_1^{-1}
 \end{array}\right\rangle.\qquad \qquad \qquad  $

\vspace{.2cm}
\noindent   Consider the third nilpotent quotient of $\pi$, \emph{i.e.} fix $k=3$. We \emph{choose} the meridians $m_i$ as (representatives of) generators, and express preferred longitudes  $\lambda_i$: we have $\lambda_3=m_ 2^{-1}n_2=m_ 2^{-1} m_1m_2m_1^{-1}$ (the other two are obtained by cyclic permutation of the indices, owing to the symmetry of the link). Taking the Magnus expansion gives 
\[ E(\lambda_3)= 1 + X_1X_2 - X_2X_1 + \textrm{terms of degree $\ge 3$}.\qquad \qquad \qquad   \]
   }
   \hspace{-1cm}
  \parbox[l]{9cm}{
  \includegraphics[scale=0.8]{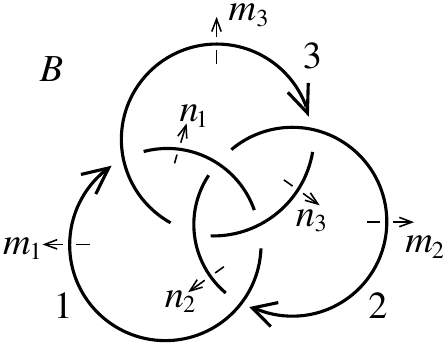}\\}
  
\noindent   We thus obtain that $\omu_B(I)=0$ for any sequences $I$ of two indices (\emph{i.e.} all linking numbers are zero), and  
$$ \omu_B(123)=\omu_B(231)=\omu_B(312)=1 \,\, \textrm{ and } \,\, \omu_B(132)=\omu_B(213)=\omu_B(321)= -1.$$
\end{example}

Example \ref{ex:bor} illustrates the fact that, just like the linking number detects (and actually, counts copies of) the Hopf link, the \emph{triple linking number} $\omu(123)$ detects the Borromean rings. This generalizes to the following realization result, due to Milnor \cite{Milnor}.

\vspace{.1cm}
\hspace{-.43cm}
\parbox[l]{10cm}{
\begin{lemma}
Let $M_n$ be the $n$-component link shown on the right \\ ($n\ge 2$), 
and let $\sigma$ be any permutation in $S_{n-1}$. Then 

\vspace{.2cm}
$\quad \omu_{M_n}\left(\sigma(1) \cdots \sigma(n-1)n(n+1)\right) = \left\{\begin{array}{cc}
1 & \textrm{if $\sigma=$Id,} \\
0 & \textrm{otherwise.}
\end{array}\right.$
\end{lemma}
}
\hspace{-5.5cm}
\parbox[l]{8cm}{
\[ \includegraphics[scale=0.95]{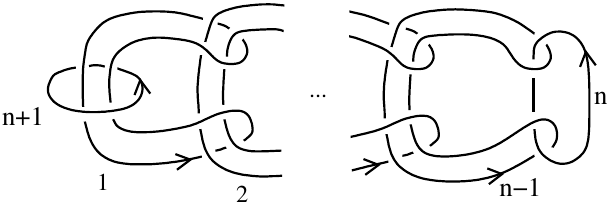} \]
}

\subsection{Some properties}
We gather here several well-known properties of Milnor $\omu$-invariants, most being due to Milnor himself \cite{Milnor,Milnor2} (unless otherwise specified). 
\subsubsection{Symmetry and Shuffle} 
It is rather clear that, if $L'$ is obtained from $L$ by reversing the orientation 
of the $i$th component, then $\omu_{L'}(I)=(-1)^{i(I)}\omu_{L}$, where $i(I)$ denotes 
the number of occurences of the index $i$ in the sequence $I$. 
The next two relations were shown by Milnor \cite{Milnor2}, 
using properties of the Magnus expansion:

 \textbf{Cyclic symmetry: }
 For any sequence $i_1\cdots i_m$ of indices in $\{1,\cdots,n\}$, we have 
  $$ \omu_L(i_1\cdots i_m)=\omu_L(i_mi_1\cdots i_{m-1}).$$
  
 \textbf{Shuffle: } 
 For any $k\in \{1,\cdots,n\}$ and any two sequences $I,J$ of indices in $\{1,\cdots,n\}$, we have 
  $$ \sum_{H\in S(I,J)} \omu_L(Hk)\equiv 0 \,\textrm{ (mod gcd$\{\Delta(Hk)$ ; $H\in S(I,J) \}$)},  $$
 where $S(I,J)$ denotes the set of all sequences obtained by inserting the indices of $J$ into $I$, preserving order. 

It follows in particular 
that there is essentially only one triple linking number: for any $3$-component algebraically split link $L$, we have  
$ \omu_L(123)=\omu_L(231)=\omu_L(312)=-\omu_L(132)=-\omu_L(213)=-\omu_L(321).$

In general, the number of independent $\omu$-invariants was given by K.~Orr, see \cite[Thm.~15]{Orr}.  

\subsubsection{Link-homotopy and concordance} \label{sec:lhc}

Milnor invariants are not only invariants of ambiant isotopy: 
they are actually invariants of \emph{isotopy}, \emph{i.e.} 
homotopy through embeddings. In particular that Milnor invariants do not see 
'local knots'. 

As mentioned in Section \ref{sec:classiflk}, the notion of link-homotopy
was introduced by Milnor himself, who proved the following in \cite{Milnor2}. 
\begin{thm}\label{thm:lh}
For any sequence $I$ of \underline{pairwise distinct} indices, $\omu(I)$ is a link-homotopy invariant. 
\end{thm}

Milnor invariants are sharp enough to detect link-homotopically trivial links, and to classify links up to link-homotopy for $\le 3$ components \cite{Milnor2}. 
The case of $4$-component links was only completed thirty years later by J.~Levine, using a refinement of Milnor's construction \cite{Levine}. The general case is discussed in Section \ref{sec:refine}. Milnor link-homotopy invariants also form a complete set of link-homotopy invariants for \emph{Brunnian} links, \emph{i.e.} links which become trivial after removal of any component (see Section \ref{sec:exmu} for examples) 
\cite{Milnor}.

Recall that two $n$-component links $L$ and $L'$ are \emph{concordant} if there is an embedding 
$ f: \sqcup_{i=1}^n (S^1\times [0,1])_i  \longrightarrow S^3\times [0,1]$ 
of $n$ disjoint copies of the annulus $S^1\times [0,1]$, such that
$f\left( (\sqcup_{i=1}^n S^1_i)\times \{ 0 \} \right)=L\times \{ 0 \}$ and $f\left( (\sqcup_{i=1}^n S^1_i)\times \{ 1 \} \right)=L'\times \{ 1 \}$. 
The following was essentially shown by J.~Stallings \cite{Stallings} (and also by A.~Casson \cite{Casson} for the more general relation of cobordism).
\begin{thm}\label{thm:concordance} 
Milnor invariants are concordance invariants. 
\end{thm}

Let us also mention here that Milnor invariants are all zero for a \emph{boundary link}, that is, for a link whose components bound mutually disjoint Seifert surfaces in $S^3$ \cite{Smythe}. 
Characterizing geometrically links with vanishing of Milnor invariants is the subject of the $k$-slice conjecture, proved by K.~Igusa and K.~Orr in \cite{IO}, which roughly states that all invariants of length $\le 2k$ vanish if and only if the link bounds a surface which ``looks like slice disks'' modulo $k$-fold commutators of the fundamental group of the surface complement. 

\subsubsection{Cabling formula} 

Milnor invariants indexed by non-repeated sequences are not only topologically relevant, by Theorem \ref{thm:lh}; they also `generate' all Milnor invariants, by the following.
\begin{thm}\label{thm:cable}
Let $I$ be a sequence of indices in $\{1,\cdots,n\}$, such that the index $i$ appears twice in $I$. 
Then $\omu_L(I) = \omu_{\tilde{L}}(\tilde{I})$,  
where $\tilde{I}$ is obtained from $I$ by replacing the second occurence of $i$ by $n+1$, and 
where $\tilde{L}$ is obtained from $L$ by adding an $(n+1)$th component, which is a parallel copy of $L_i$ having linking number zero with it.
\end{thm}

\begin{example}\label{ex:W}
A typical example of a link-homotopically trivial link is the Whitehead link $W$, shown below 

\hspace{-6.2cm}
\parbox[l]{3cm}{
\[ \includegraphics[scale=0.9]{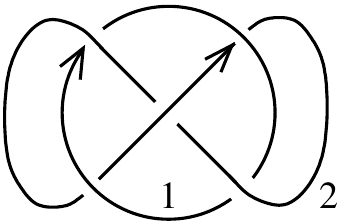} \,\,\,\,  \]
}
\parbox[l]{8cm}{ 
on the left-hand side. 
But Milnor concordance invariants do detect $W$; specifically, we have $\omu_W(1122)=1$. 

As a matter of fact, one  can check that 
$$\hspace{-6.5cm} \omu_W(1122)= \omu_{\tilde{W}}(1324)=1,$$
where $\tilde{W}$  is the $4$-component link shown on the right. 
}
\hspace{-5.9cm}
\parbox[l]{3cm}{
\[ \includegraphics[scale=0.8]{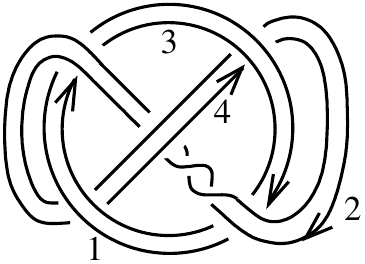} \]
}
\end{example}

\subsubsection{Additivity} 
Given two $n$-component links $L$ and $L'$, lying in two disjoint $3$-balls of $S^3$, a \emph{band sum} 
 $L\sharp L'$ of $L$ an $L'$ is any link made of pairwise disjoint connected sums of the $i$th components of $L$ and $L'$ ($1\le i\le n$). 
Although this operation is not well-defined, V.~Krushkal showed that, for any sequence $I$, 
 $$ \omu_{L\sharp L'}(I) \equiv \omu_{L}(I) + \omu_{L'}(I)  \,\left( \textrm{mod gcd$(\Delta_L(I),\Delta_{L'}(I))$}\right), $$
meaning that Milnor invariants are independent of the choice of bands defining $L\sharp L'$, see  \cite{Krushkal}.

\section{Further properties}\label{sec:more}

Since their introduction in the fifties, Milnor invariants have been the subject of numerous works. 
In this section, which claims neither exhaustivity nor precision, we briefly overview some of these results.

\subsection{Refinements and generalizations}\label{sec:refine}
The main difficulty in understanding Milnor $\omu$-invariants lies in these intricate indeterminacies $\Delta(I)$. Multiple attempts have been made to refine this indeterminacy, in order to get more subtle invariants. 

T.~Cochran defined, by considering recursive intersection curves of Seifert surfaces, link invariants which recover  (and shed beautiful geometric lights on) $\omu$-invariants, and sometimes refine them \cite{Cochran_AMS}. 

K.~Orr defined invariants, as the class of the ambient $S^3$ in some third homotopy group, which refines significantly the indeterminacy of Milnor invariants \cite{Orr1,Orr}. Orr's invariant is also the first attempt towards a transfinite version of Milnor invariants, a problem posed in \cite{Milnor2}; see also the work of J.~Levine in \cite{Levine_transfinite}. 

A decisive step was taken by Habegger and Lin, who showed that the indeterminacies $\Delta_L(I)$
are equivalent to the indeterminacy in representing $L$ as the closure of a \emph{string link}, \emph{i.e.} a pure tangle without closed component \cite{HL} (see also \cite{levine_based}). 
This led them to a full link-homotopy classification of (string) links \cite{HL}. 

\subsection{Link maps and higher dimensional links} 
There are several higher dimensional versions of the linking number; see \cite[\S~3.4]{CKS} for a good survey.
Some are invariants of \emph{link maps}, \emph{i.e.} maps from a union of two  spheres (of various dimensions) to a sphere with disjoint images 
-- these are the natural objects to consider when working up to link-homotopy. 
In particular, the first example of a link map which is not link-homotopically trivial was given in 
\cite{FR} using an appropriate generalization of the linking number. \\
Higher dimensional generalizations of Milnor invariants were defined and extensively studied by U.~Koschor\-ke for link maps with many components, 
see \cite{Ko,KoLink} and references therein. 
Note also that Orr's invariants, which generalize Milnor invariants (see Section \ref{sec:refine}), are also defined in any dimensions.  
But there does not seem to be nontrivial analogues of Milnor invariants for $n$-dimensional links ($n\ge 2$), \emph{i.e.} for embedded $n$-spheres in codimension $2$; in fact, all are link-homotopically trivial \cite{BT}. 
However, Milnor invariants generalize naturally to $2$-string links, 
\emph{i.e.} knotted annuli in the $4$-ball bounded by a prescribed unlink in the $3$-sphere, 
and classify them up to link-homotopy \cite{AMW}, thus providing a higher dimensional version of \cite{HL}.

\subsection{Relations with other invariants}

Milnor invariants are not only natural generalizations of the linking number, but are also directly related to this invariant. Indeed, K.~Murasugi expressed Milnor $\omu$-invariants of a link as a linking number in certain branched coverings of $S^3$ along this link \cite{Mu85}. 
The Alexander polynomial is also rather close in nature to Milnor invariants, being extracted from the fundamental group of the complement. 
As a matter of fact, there are a number of results relating these two invariants: see for example  \cite{C,L,Mu,Polyak,T}. 
On the other hand, the relation to the Konstevich integral \cite{HM} hints to potential connections to quantum invariants. 
Such relations were given with the HOMFLY-PT polynomial \cite{MYGT} and the (colored) Jones polynomial \cite{MSu}. 
Milnor string link invariants also satisfy a skein relation \cite{Polyak_skein}, which is a typical feature of polynomial and quantum invariants. 

There also are known relations outside knot theory. 
Milnor $\omu$-invariants can be expressed in terms of Massey products of the complement, which are higher order cohomological invariants generalizing the cup product \cite{Porter,Turaev}.
Also, Milnor string link invariants are in natural correspondence with Johnson homomorphisms of homology cylinders, which are $3$-dimensional extensions of certain abelian quotients of 
the mapping class group \cite{Habegger}.
Finaly, as part of the deep analogies he established between knots and primes,  
M.~Morishita defined and studied arithmetic analogues of Milnor invariants for prime numbers \cite{Morishita}. 

\subsection{Finite type invariants} 
The linking number is (up to scalar) the unique degree $1$ finite type link invariant. 
This was generalized independently by D.~Bar Natan \cite{BN_SL} and X.S.~Lin \cite{Lin}, who showed that any Milnor \emph{string link} invariant of length $k$ is a finite type invariant of degree $k-1$. 
As a consequence, Milnor string link invariants can be, at least in principle, extracted from the Kontsevich integral, which is universal among finite type invariants. 
This was made completely explicit by G.~Masbaum and N.~Habegger in \cite{HM}. 

Note that the finite type property does not make sense for higher order $\omu$-invariants (of \emph{links}), since the indeterminacy $\Delta(I)$ is in general not the same for two links which differ by a crossing change. 
Nonetheless,  $\omu$-invariants of length $k$ are invariants of $C_{k}$-equivalence \cite{Habiro}, a property shared by all degree $k-1$ invariants. 

\subsection{Virtual theory}

There are two distinct extensions of the linking number for virtual $2$-component links, namely $\lk_{1,2}$ and $\lk_{2,1}$,  where $\lk_{i,j}$ is the sum of signs of crossings where $i$ passes over $j$. 
Notice that these virtual linking numbers are actually invariants of \emph{welded links}, \emph{i.e.} are invariant under the forbidden move allowing a strand to pass \emph{over} a virtual crossing. 
Indeed, the linking number is extracted from (a quotient of) the fundamental group, which is a welded invariant \cite{Kauffman}. Likewise, extensions of Milnor invariants shall be welded invariants. 
A general welded extension of Milnor string link invariants is given in \cite{ABMW}, which classifies welded string links up to self-virtualization,  generalizing the classification of \cite{HL}. 
This extension recovers and extends that of \cite{PK}, which gives general Gauss diagram formulas for virtual Milnor invariants. 

\begin{acknowledgments}
The author would like to thank Benjamin Audoux and Akira Yasuhara, for their comments on a preliminary version of this note. 
\end{acknowledgments}
\bibliographystyle{abbrv}
\bibliography{lkmilnor}

\end{document}